\documentclass[a4paper]{scrartcl}

\usepackage[utf8]{inputenc}
\usepackage[T1]{fontenc}
\usepackage{lmodern}

\usepackage{amsmath}
\usepackage{amssymb}
\usepackage{amsthm}

\usepackage{nicefrac}

\theoremstyle{plain}
\newtheorem{theorem}{Theorem}
\newtheorem{lemma}[theorem]{Lemma}
\newtheorem{corollary}[theorem]{Corollary}

\theoremstyle{definition}
\newtheorem{remark}[theorem]{Remark}

\usepackage{tikz}
\usepackage[hypertexnames=false]{hyperref}

\DeclareMathOperator{\st}{\,s.t.}

\newcommand{\sigp}[1]{\Sigma^{\text P}_{#1}}

\newcommand{\R}{{\mathbb{R}}}
\newcommand{\U}{{\mathcal{U}}}

\begin{document}

\title{On the Complexity of Robust Bilevel Optimization With Uncertain Follower's Objective\thanks{This work has partially been supported by Deutsche Forschungsgemeinschaft (DFG) under grant no.~BU 2313/6.}}

\author{Christoph Buchheim, Dorothee Henke, Felix Hommelsheim\\[1ex] Department of Mathematics, TU Dortmund University, Germany}

\date{}

\maketitle

\begin{abstract}
  We investigate the complexity of bilevel combinatorial
  optimization with uncertainty in the follower's objective, in a
  robust optimization approach. We show that the robust counterpart
  of the bilevel problem under interval uncertainty can
  be~$\sigp 2$-hard, even when the certain bilevel problem is
  NP-equivalent and the follower's problem is tractable. On the
  contrary, in the discrete uncertainty case, the robust bilevel
  problem is at most one level harder than the follower's problem.
 
  \smallskip
  
\end{abstract}

\section{Introduction}

When addressing combinatorial optimization problems, it is usually
assumed that there is only one decision maker that has control over
all variables, and that all parameters of the problem are known
precisely.  However, these two assumptions are not always satisfied in
real-world applications. The latter issue of data uncertainty has led
to the growing research fields of \emph{stochastic} and \emph{robust
  optimization}. In this paper, we focus on the robust optimization
approach to combinatorial optimization and assume that only objective
functions are uncertain. The complexity of the resulting robust
counterparts highly depends on the chosen \emph{uncertainty sets},
which are supposed to contain all possible scenarios, while the
objective function value in the worst-case scenario is optimized.
Typical classes of uncertainty sets are finite sets \emph{(discrete
  uncertainty)} or \mbox{hyper} boxes \emph{(interval uncertainty)}. In the
standard setting, interval uncertainty preserves the complexity of the
problem, while discrete uncertainty often leads to NP-hard robust
counterparts even when the underlying certain problem is
tractable. For these and other complexity results, we refer
to~\cite{kouvelis} or to the recent survey~\cite{bk17}. For a general
introduction to robust optimization, see, e.g., \cite{ben2002robust}.

A common approach to address optimization problems involving more than
one decision maker is \emph{multilevel optimization} (or \emph{bilevel
  optimization} in case of two decision makers).  A multilevel
optimization problem models the interplay between several decision
makers, each of them having their own decision variables, objective
function and constraints. The decisions may depend on each other and
are made in a hier\-ar\-chi\-cal order: first, the highest-level
decision maker decides.  Based on this, the decision maker on the next
level takes his choice, followed by the third-level decision maker,
and so on.  Usually, the problem is viewed from the perspective of the
highest-level decision maker, who has perfect knowledge about all the
lower-level problems.  In particular, she has to anticipate the entire
sequence of optimal lower-level responses to her decision.  One has to
distinguish between two kinds of behavior in case some lower-level
decision maker has more than one optimal solution: in the
\emph{optimistic} setting, he acts in favor of the highest-level
decision maker (or another upper-level decision maker), and in the
\emph{pessimistic} setting, he chooses a worst-case solution for her.

Jeroslow~\cite{jeroslow85} showed that an $\ell$-level program with
only linear constraints and objectives is already $\sigp
{\ell-1}$-hard in general, and thus NP-hard in the bilevel case, where
he assumes the optimistic setting for every pair of consecutive
players.  Furthermore, if the variables of all decision makers are
restricted to be binary, then the problem turns out to be $\sigp
{\ell}$-hard in general.  Many bilevel variants of classical
combinatorial optimization problems have also shown to be NP-hard,
e.g., a bilevel assignment problem~\cite{gassner2009computational} or
a bilevel minimum spanning tree problem~\cite{buchheim2020complexity}.
For a thorough overview of bilevel and multilevel optimization, we
refer to the surveys~\cite{vicente94,colson07,dempe03,lu16} or to the
book~\cite{dempe15}.

In this paper, our aim is to bring together robust and bilevel
optimization and to investigate the result in terms of
complexity. More precisely, we address bilevel combinatorial
optimization problems where the lower-level objective function is
unknown to the upper-level decision maker.  However, the latter knows
an uncertainty set to which this objective function will belong, and
she aims at optimizing her worst-case objective.  While the two
decision makers are often referred to as leader and follower,
motivated by the origins of bilevel optimization in Stackelberg
games~\cite{Stackelberg}, our situation can be illustrated by
considering three decision makers: the \emph{leader} first chooses a
solution, then an \emph{adversary} chooses an objective function for
the follower, and the \emph{follower} finally computes his optimal
solution according to this objective function and depending on the
leader's choice. The aim of the leader is to optimize her own
objective value, which depends on the follower's choice, while the
adversary has the opposite objective.

We start our investigation with a certain combinatorial bilevel
problem of the form
\begin{equation}\label{eq:basic_problem}\tag{P}
\begin{aligned}
  \max_x~ & d^\top y_c\\
  \st~ &  x\in X \\
  & y_c \in
  \begin{aligned}[t]
    \arg\max_y~ & c^\top y\\
    \st~ & y\in Y(x)
  \end{aligned}
\end{aligned}
\end{equation}
with~$X\subseteq\{0,1\}^p$, $Y(x)\subseteq\R^n$ for all~$x \in X$, and
$c, d \in \R^n$. We note that the notation of~\eqref{eq:basic_problem}
(and all other bilevel programs throughout this paper) is convenient,
but not precise in case the follower does not have a unique optimum
solution. Formally, the optimistic or pessimistic case could be
modeled by adding a maximization or minimization over $y_c$,
respectively, to the leader's objective function.

We concentrate on the most basic setting where the follower's problem
is a linear program with only the right hand side depending
(affine-linearly) on the leader's decision~$x$. More formally, we
assume that there are~$A\in\R^{m\times n}$, $B\in\R^{m\times p}$,
and~$b\in\R^m$ such that~$Y(x)=\{y\in\R^n\mid Ay\le Bx+b\}$ for all~$x
\in X$.  For simplicity, we assume throughout that~$Y(x)$ is non-empty
and bounded for all~$x\in X$, so that an optimal follower's response
always exists.

In particular, for given~$x$ and~$c$, the follower's problem is a
linear program and thus always tractable. This holds not only in terms
of finding \emph{any} optimal solution, but also if the task includes
finding an optimal solution that is best or worst possible for the
leader. This amounts to optimizing two linear objective functions
lexicographically, which can be done
efficiently~\cite{SheraliSoyster}.  It follows that the leader's
optimization problem~\eqref{eq:basic_problem} is NP-easy, since its
decision variant belongs to NP, the solution~$x\in X$ being the
certificate.

Turning our attention to the robust optimization setting, we now
assume that the follower's objective~$c$ is uncertain. However, we are
given an uncertainty set~$\U\subseteq\R^n$ containing all relevant
realizations of~$c$. The robust counterpart
of~\eqref{eq:basic_problem} then reads
\begin{equation}\label{eq:robust_problem}\tag{R}
\begin{aligned}
  \max_x~ & \min_{c\in\U}~d^\top y_c\\
  \st~ &  x\in X \\
  & y_c \in
  \begin{aligned}[t]
    \arg\max_y~ & c^\top y\\
    \st~ & y\in Y(x)\;.
  \end{aligned}
\end{aligned}
\end{equation}
Our aim is to investigate the complexity of~\eqref{eq:robust_problem}
relative to the complexity of~\eqref{eq:basic_problem}. In case the
follower's solution~$y_c$ for a given objective~$c$ is not unique, we
will consider both the optimistic and the pessimistic approach.

Besides the robust counterpart~\eqref{eq:robust_problem}, we will also
consider the adversary's subproblem, i.e., the problem of evaluating
the leader's objective function in~\eqref{eq:robust_problem}. For a
fixed leader's choice~$x\in X$, this problem thus reads
\begin{equation}\label{eq:adversary_problem}\tag{A}
\begin{aligned}
  \min_c~ & d^\top y_c\\
  \st~ &  c\in\U \\
  & y_c \in
  \begin{aligned}[t]
    \arg\max_y~ & c^\top y\\
    \st~ & y\in Y(x)\;,
  \end{aligned}
\end{aligned}
\end{equation}
which turns out to be a bilevel optimization problem again, however
with the leader's constraints having a very specific structure now,
depending only on the uncertainty set~$\U$.  Finally, for fixed~$x$
and~$c$, we will consider the certain follower's problem
\begin{equation}\label{eq:follower_problem}\tag{F}
\begin{aligned}
  \max_y~ & c^\top y\\
  \st~ & y\in Y(x)\;,
\end{aligned}
\end{equation}
which, as already mentioned, is a linear program in our setting.

It is easily verified that the robust
problem~\eqref{eq:robust_problem} can be at most one level harder, in
the polynomial hierarchy, than the adversary's
problem~\eqref{eq:adversary_problem}, i.e., the evaluation of the
objective function of~\eqref{eq:robust_problem}. The complexity
of~\eqref{eq:robust_problem} and~\eqref{eq:adversary_problem} strongly
depends on the uncertainty set~$\U$.  We will see
that~\eqref{eq:robust_problem} can be $\sigp 2$-hard
and~\eqref{eq:adversary_problem} NP-hard in case of interval
uncertainty, while in case of discrete uncertainty,
\eqref{eq:adversary_problem} is tractable and thus
\eqref{eq:robust_problem} is NP-easy.  This is a remarkable difference
to classical (single-level) robust optimization, where, as mentioned
above, interval uncertainty in the objective function does not make a
problem harder, while many tractable combinatorial optimization
problems have an NP-hard robust counterpart when the set~$\U$ is
finite.

Addressing similar complexity questions, Buchheim and
Henke~\cite{knapsack} investigate a bilevel continuous knapsack
problem where the leader controls the continuous capacity of the
knapsack and, as in~\eqref{eq:robust_problem}, the follower's profits
are uncertain. They show that (the analoga of)
both~\eqref{eq:follower_problem} and~\eqref{eq:basic_problem} can be
solved in polynomial time in this case,
while~\eqref{eq:adversary_problem} and~\eqref{eq:robust_problem} turn
out to be NP-hard for, e.g., budgeted and ellipsoidal uncertainty. On
the other hand, both \eqref{eq:adversary_problem}
and~\eqref{eq:robust_problem} with respect to discrete and interval
uncertainty remain tractable. The latter result rises the question
whether or not interval or discrete uncertainty, in general, can
increase the complexity when going from~\eqref{eq:basic_problem}
to~\eqref{eq:robust_problem}.

The complexity of bilevel optimization problems under uncertainty, in
a robust optimization framework, has also been addressed in a few
recent articles. However, the assumptions about what is uncertain are
different to our setting. In~\cite{Besancon2020}, multilevel
optimization problems are investigated in which some follower's
decision does not have to be optimal according to his objective, but
can deviate from the optimum value by a small amount; see
also~\cite{Besancon2019} and the references
therein. In~\cite{Beck2021}, bilinear bilevel programming problems are
considered and it is assumed that the follower cannot observe the
leader's decision precisely. In both settings, it turns out that the
complexity does not increase significantly with respect to the
corresponding model without uncertainty.

\section{Interval uncertainty}

We start by considering the case of interval uncertainty, i.e.,
uncertainty sets of the form $\U = [c_1^-, c_1^+] \times \dots \times
[c_n^-, c_n^+]$ for given $c^-, c^+ \in \R^n$ with $c^- \le c^+$. We
will prove that, in this case, the robust
counterpart~\eqref{eq:robust_problem} is significantly harder than the
underlying problem~\eqref{eq:basic_problem} in general, while the
adversary's problem~\eqref{eq:adversary_problem} is significantly
harder than the follower's problem~\eqref{eq:follower_problem}.  For
the reductions, we will use the well-known \emph{satisfiability
  problem} (SAT) as well as the less-known \emph{quantified
  satisfiability problem} (QSAT), defined as follows:

\medskip

(SAT) Given a Boolean formula $f \colon \{0, 1\}^n \to \{0, 1\}$,
decide whether there exists a satisfying assignment, i.e., a~$y \in
\{0, 1\}^n$ such that $f(y) = 1$.

\medskip

(QSAT) Given a Boolean formula $f \colon \{0, 1\}^p \times \{0, 1\}^n
\to \{0, 1\}$, decide whether there exists some~$x \in \{0, 1\}^p$
such that $f(x, y) = 1$ for all $y \in \{0, 1\}^n$.

\medskip

For both problems, we may assume that the function~$f$ is given in the
input by a recursive definition using negations and bivariate
conjunctions and disjunctions. The problem SAT is well-known to be
NP-complete, while QSAT is an important example of a $\sigp
2$-complete problem~\cite{stockmeyer77}. The complexity class~$\sigp
2$ contains all decision problems which can be solved by a
nondeterministic Turing machine equipped with an oracle for some
NP-complete decision problem; see~\cite{stockmeyer77} for an
introduction to the polynomial-time hierarchy and a formal definition
of the complexity classes $\sigp \ell$.  For later use, we also recall
the following definitions: a decision or optimization problem is
\emph{NP-easy} if it can be polynomially reduced to an NP-complete
decision problem, and \emph{NP-equivalent} if it is both NP-hard and
NP-easy.

\begin{theorem}\label{theorem:qsat}
  In case of interval uncertainty, the robust
  counterpart~\eqref{eq:robust_problem} can be $\sigp 2$-hard and the
  adversary's problem~\eqref{eq:adversary_problem} NP-hard, even
  if~$X=\{0,1\}^p$.
\end{theorem}

\begin{proof}
  We prove the $\sigp 2$-hardness of~\eqref{eq:robust_problem} by a
  polynomial reduction from the $\sigp 2$-complete problem
  QSAT. Since~\eqref{eq:robust_problem} can be at most one level
  harder than~\eqref{eq:adversary_problem}, as mentioned in the
  introduction, the NP-hardness of~\eqref{eq:adversary_problem}
  directly follows. Furthermore, a reduction from the complement of
  SAT to~\eqref{eq:adversary_problem} is implicitly contained in our
  reduction from QSAT to~\eqref{eq:robust_problem}.

  For a given Boolean formula~$f\colon \{0,1\}^p \times \{0,
  1\}^n\to\{0,1\}$, consider the bilevel problem
  \begin{equation}\label{eq:qsat}\tag{Q}
    \begin{aligned}
      \max_x~ & \min_{c\in [-1,1]^n}~f(x,y_c)\\
      \st~ &  x\in \{0,1\}^p \\
      & y_c \in
      \begin{aligned}[t]
        \arg\max_y~ & c^\top y\\
        \st~ & y\in [0,1]^n\;.
      \end{aligned}
    \end{aligned}
  \end{equation}
  This problem can be linearized in the standard way: introduce a new
  variable with domain~$[0,1]$ for each intermediate term in the
  recursive definition of~$f$, including~$f$ itself, and model the
  corresponding nonlinear relations between all variables by
  appropriate linear inequalities (for conjunctions and disjunctions)
  or linear equations (for negations); see, e.g.,
  \cite{jeroslow85}. In order to obtain a reformulation
  of~\eqref{eq:qsat} in the form of~\eqref{eq:robust_problem}
  with~$X=\{0,1\}^p$, we add all new variables and linear constraints
  to the follower's problem.  For all newly introduced variables, we
  set the corresponding entries of~$c^-$ and~$c^+$ to $0$.

  We first consider the optimistic setting.  Assume that~$f$ is a
  no-instance of~QSAT. Then for all leader's choices~$x\in\{0,1\}^p$,
  there exists some~$y\in\{0,1\}^n$ such that~$f(x,y)=0$. By
  construction of~\eqref{eq:qsat}, the adversary can enforce any
  follower's response~$y\in\{0,1\}^n$ by appropriately
  choosing~$c\in\{-1,1\}^n$. In that case, since~$x$ and~$y$ are
  binary, the linearization forces all additional variables to the
  corresponding binary values as well, including the variable
  representing the objective function. This shows that the optimal
  value of~\eqref{eq:qsat} is $0$ in this case.

  If~$f$ is a yes-instance, there exists~$x\in\{0,1\}^p$ such
  that~$f(x,y)=1$ for all~$y\in\{0,1\}^n$. Then we claim that this
  leader's choice~$x$ yields an objective value of~$1$ in the
  linearization of~\eqref{eq:qsat}. Indeed, for every adversary's
  choice~$c$, there is an optimal follower's solution such that~$y_c$
  is binary. By construction, each such solution leads
  to~$f(x,y_c)=1$. Since we adopt the optimistic approach, the
  follower will thus choose some solution~$y_c$ with~$f(x,y_c)=1$.

  Finally, assume that we adopt the pessimistic approach. In this
  case, we need to modify the construction of~\eqref{eq:qsat} as
  follows. For every follower's variable~$y_i$, we introduce a new
  follower's variable~$\bar y_i$, along with linear constraints
  $$\bar y_i\ge 0,~\bar y_i\le y_i,~\bar y_i\le 1-y_i\;.$$ Each new
  variable~$\bar y_i$ has coefficient~$M$ in the leader's objective,
  where $M \ge 3$ is at least the number of atomic terms in $f$, and a
  certain coefficient of~$1$ in the follower's objective.  This
  ensures that $\bar y_i = \min\{y_i, 1-y_i\}$ holds in every optimal
  follower's solution, i.e., the variable~$\bar y_i$ models the
  deviation of $y_i$ from the closest binary value.

  \begin{figure}
    \centering
    \begin{tikzpicture}[scale=0.85]
      \draw[->] (0,0) -- (3.5,0) node[below] {$y_i$};
      \draw[->] (0,-2.2) -- (0,2.5) node[left] {$c_i y_i + \bar y_i$};
      \draw (1,0.1) -- (1,-0.1) node[below] {$\frac 12$};
      \draw (2,0.1) -- (2,-0.1) node[below] {$1$};
      \draw (0.1,0) -- (-0.1,0) node[left] {$0$};
      \draw (0.1,1) -- (-0.1,1) node[left] {$\frac 12$};
      \draw (0.1,2) -- (-0.1,2) node[left] {$1$};
      \draw (0.1,-1) -- (-0.1,-1) node[left] {$- \frac 12$};
      \draw (0.1,-2) -- (-0.1,-2) node[left] {$-1$};

      \draw[very thick] (0,0) -- (1,2) -- (2,2) node[right] {\small $c_i = 1$};
      \draw[very thick] (0,0) -- (1,1.5) -- (2,1) node[right] {\small $c_i = \frac 12$};
      \draw[very thick] (0,0) -- (1,1) -- (2,0) node[above right] {\small $c_i = 0$};
      \draw[very thick] (0,0) -- (1,0.5) -- (2,-1) node[right] {\small $c_i = - \frac 12$};
      \draw[very thick] (0,0) -- (1,0) -- (2,-2) node[right] {\small $c_i = -1$};
    \end{tikzpicture}
    \caption{Illustration of the follower's objective function in the
      pessimistic setting in the proof of
      Theorem~\ref{theorem:qsat}. For a fixed dimension $i \in \{1,
      \dots, n\}$, the value $c_i y_i + \bar y_i$ is displayed in
      dependence of the value $y_i \in [0,1]$ the follower chooses,
      for different values of $c_i$ chosen by the adversary.}
    \label{fig:follower_objective_pessimistic}
  \end{figure}
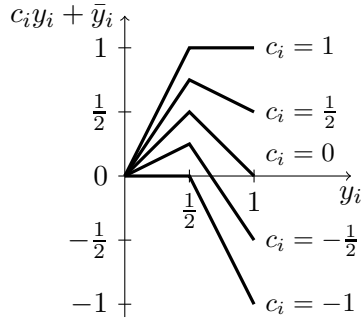

  The resulting follower's objective $c_i y_i + \bar y_i$, in every
  dimension $i \in \{1, \dots, n\}$, depends on the adversary's choice
  $c_i$ and always consists of two linear pieces, as a function
  in~$y_i$; see Figure~\ref{fig:follower_objective_pessimistic}.  If
  the adversary chooses~$c\in\{-1,1\}^n$, the follower's optimum
  solutions satisfy $y_i = \bar y_i \in [0, \nicefrac 12]$ or $y_i = 1
  - \bar y_i \in [\nicefrac 12, 1]$, respectively.  By setting all
  $y_i$ to binary values, the follower can achieve a leader's
  objective value of at most $1$ then, due to the pessimistic
  assumption.  If the adversary chooses some~$c_i\in(-1,1)$, then by
  construction, the value of~$(y_i,\bar y_i)$ in any optimal
  follower's solution is~$(\nicefrac 12,\nicefrac 12)$, so that the
  leader's objective is at least~$\nicefrac M2> 1$, which is never
  optimal for the adversary.  Hence, the adversary must always
  choose~$c\in\{-1,1\}^n$.

  For a no-instance~$f$ and any~$x\in\{0,1\}^p$, the adversary can now
  choose $c \in \{-1, 1\}^n$ such that the follower has a binary
  vector~$y \in \{0, 1\}^n$ in his set of optimum solutions for which
  $f(x, y) = 0$ holds. By the pessimistic assumption, the follower
  actually chooses $y$ such that the leader's objective value is $0$.

  Assume, on the contrary, that~$f$ is a yes-instance, and consider
  some leader's solution~$x\in\{0,1\}^p$ such that~$f(x,y)=1$ for
  all~$y\in\{0,1\}^n$.  Any adversary's choice $c \in \{-1, 1\}^n$,
  together with any binary follower's choice $y$ then results in an
  objective value of $1$ for the leader.  Since we adopt the
  pessimistic approach, it remains to show that the follower cannot
  achieve a smaller objective value for the leader by choosing $y_i$
  nonbinary and exploiting the resulting flexibility when setting the
  linearization variables. However, one can verify that this is not
  possible, because the gain in the variable representing~$f(x, y)$
  would be exceeded by the punishment due to the terms $M \bar y_i$,
  by definition of~$M$; see, e.g., Lemma~4.1
  in~\cite{jeroslow85}. Thus, the follower chooses $y \in \{0, 1\}^n$
  and the leader's optimal value is indeed $1$.
\end{proof}

Note that the underlying certain problem of~\eqref{eq:qsat} is
NP-equivalent. Indeed, it is NP-hard even when~$n=0$, i.e., when the
follower only controls the additional variables introduced for
linearizing~$f$, which are uniquely determined by the leader's
choice~$x\in\{0,1\}^p$. The problem is then equivalent to~SAT. On the
other hand, the problem is NP-easy, for any $n$, as argued in the
introduction. We emphasize however that the optimal follower's choice
as well as the leader's optimal value in the underlying certain
problem might not be binary if~$n\neq 0$, for some values of $c$.

In order to strengthen the statement of Theorem~\ref{theorem:qsat},
the integrality constraints in the leader's problem can be relaxed
without changing the result, using a similar construction as in the
pessimistic case (and in~\cite{jeroslow85}): for each leader's
variable~$x_i$, we can introduce a new follower's variable~$\bar
x_i\in[0,1]$ with coefficient~$-M$ in the leader's objective,
where~$M\ge3$ is again at least the number of atomic terms in~$f$,
with certain coefficient~$1$ in the follower's objective, and with
additional follower's constraints $\bar x_i \ge 0,~\bar x_i\le
x_i,~\bar x_i\le 1-x_i$. This ensures~$\bar x_i=\min\{x_i,1-x_i\}$ in
any optimal follower's solution and, by the choice of $M$,
that~$x\in\{0,1\}^p$ holds in any optimal leader's solution.  We
obtain a problem of type~\eqref{eq:robust_problem} again, but without
integrality constraints. Note that this construction does not
interfere with the one in the proof of Theorem~\ref{theorem:qsat} in
the pessimistic case because the latter is relevant only in case of a
yes-instance, whereas the former is relevant only in case of a
no-instance.

\begin{remark}
  If binarity constraints were allowed also in the follower's problem,
  a similar construction as in Theorem~\ref{theorem:qsat} would
  actually yield a $\sigp 3$-hard problem, in the optimistic
  setting. For this, consider the $\sigp 3$-complete problem

  \medskip
  
  (QSAT$_3$) Given a Boolean formula $f \colon \{0, 1\}^p \times \{0,
  1\}^n \times \{0, 1\}^m \to \{0, 1\}$, decide whether there exists
  some~$x \in \{0, 1\}^p$ such that for all $y \in \{0, 1\}^n$ there
  exists some $z\in \{0, 1\}^m$ with $f(x, y, z) = 1$.

  \medskip

  \noindent One could reduce QSAT$_3$ to~\eqref{eq:robust_problem} by
  considering two separate sets of variables in the follower's
  problem, the first one corresponding to~$y$ and being ``controlled''
  by the adversary (i.e., with~$c_i^-=-1$, $c_i^+=1$) and the second
  one corresponding to~$z$ and being ``controlled'' by the follower
  (i.e., with~$c_i^-=c_i^+=0$). The adversary's problem is then
  equivalent to the complement of QSAT and thus~$\sigp 2$-hard. The
  underlying certain problem is equivalent to SAT, similarly as in
  Theorem~\ref{theorem:qsat}, and thus NP-equivalent. However, in
  contrast to the setting of Theorem~\ref{theorem:qsat}, the same
  holds for the follower's problem here. Hence, the increase in
  complexity of the robust problem relative to the underlying certain
  problem is strictly greater here than in Theorem~\ref{theorem:qsat},
  while it remains the same relative to the follower's problem.
\end{remark}

\begin{remark}\label{rem:int}  
  The result of Theorem~\ref{theorem:qsat} does not only hold for
  interval uncertainty sets, but also for many other uncertainty
  sets~$\U$ where each entry of~$c$ can be chosen to be positive or
  negative independently of each other. In particular, the proof of
  Theorem~\ref{theorem:qsat} can be easily adapted for the case of
  uncorrelated discrete uncertainty, as considered in~\cite{knapsack}.
\end{remark}

\section{Discrete Uncertainty}

We next address the case of a finite uncertainty set~$\U$. For this
case, we first show that robust counterparts of bilevel problems are
at least as hard as robust counterparts of single-level problems, for
each underlying feasible set of the leader.
\begin{lemma}\label{lemma:discrete}
  For any class of underlying sets~$X\subseteq\{0,1\}^p$, the robust
  counterpart of linear optimization over~$X$ subject to discrete
  uncertainty can be polynomially reduced to a problem of
  type~\eqref{eq:robust_problem}, with the same number of scenarios.
\end{lemma}
\begin{proof}
  Starting from the robust single-level counterpart
  \begin{equation}\label{eq:robust_problem_single}\tag{S}
    \begin{aligned}
      \max_x~ & \min_{c\in\U}~c^\top x\\
      \st~ &  x\in X
    \end{aligned}
  \end{equation}
  of linear optimization over some set $X \subseteq \{0,1\}^p$,
  with~$\U=\{c_1,\dots,c_m\} \subseteq \R^p$, we construct the robust
  bilevel problem
  \begin{equation}\label{eq:robust_problem_red}\tag{R(S)}
    \begin{aligned}
      \max_x~ & \min_{\tilde c\in \tilde\U}~y_{\tilde c}\\
      \st~ &  x\in X \\
      & (y_{\tilde c},z_{\tilde c}) \in
      \begin{aligned}[t]
        \arg\max_{y, z}~ & \tilde c^\top z\\[-1ex]
        \st~ & y=\sum_{j=1}^m \sum_{i=1}^p c_{ji} x_i z_j\\[-1ex]
        & y\in\R,\quad z\in \R^m_{\ge0},\quad \sum_{j=1}^m z_j=1
      \end{aligned}
    \end{aligned}
  \end{equation}
  with uncertainty set~$\tilde\U=\{e_1,\dots,e_m\} \subseteq \R^m$,
  where $e_i$ denotes the $i$-th unit vector.  When the adversary
  chooses~$\tilde c=e_j$, the follower's unique optimum solution
  consists of~$z=e_j$ and~$y=c_j^\top y$.
  Thus~\eqref{eq:robust_problem_single}
  and~\eqref{eq:robust_problem_red} are equivalent.

  It remains to show that we can
  linearize~\eqref{eq:robust_problem_red} and thus turn it into the
  form required in~\eqref{eq:robust_problem}.  For this, we introduce
  follower's variables $u_{ji}$ and constraints
  $$u_{ji}\ge 0,~u_{ji}\ge x_i+z_j-1,~u_{ji}\le x_i,~u_{ji}\le z_j$$
  for all~$j\in\{1,\dots,m\}$ and~$i\in\{1,\dots,p\}$, and replace the
  first constraint in the follower's problem by
  $$y=\sum_{j=1}^m \sum_{i=1}^p c_{ji} u_{ji}\;.$$
  For every~$\tilde c \in \tilde\U$, the unique optimum follower's
  solution satisfies~$z\in\{0,1\}^m$.  Hence this linearization
  ensures~$u_{ji}=x_iz_j$ and thus yields a problem equivalent
  to~\eqref{eq:robust_problem_red}, having the form
  of~\eqref{eq:robust_problem}.  Formally, all scenarios~$\tilde
  c\in\tilde \U$ must be extended by zeros in the entries
  corresponding to the variables~$y$ and~$u_{ji}$.
\end{proof}
Note that all optimal follower's solutions in the proof of
Lemma~\ref{lemma:discrete} are uniquely determined, so that the result
holds for both the optimistic and the pessimistic setting.

The robust counterpart of a single-level combinatorial optimization
problem under discrete uncertainty often turns out to be NP-hard even
for two scenarios, and strongly NP-hard when the number of scenarios
is not fixed; this holds true even for unrestricted binary
optimization~\cite{bk17}. Together with Lemma~\ref{lemma:discrete},
this immediately implies
\begin{corollary} \label{cor:discrete}
  Problem~\eqref{eq:robust_problem} with $X=\{0,1\}^p$ can be NP-hard
  for~$|\U|=2$ and strongly NP-hard for finite~$\U$, even
  when~\eqref{eq:basic_problem} is solvable in polynomial time.
\end{corollary}

On the other hand, the complexity of~\eqref{eq:robust_problem} can be
bounded from above as follows: as mentioned in the
introduction,~\eqref{eq:robust_problem} can be at most one level
harder
than~\eqref{eq:adversary_problem}. Moreover,~\eqref{eq:adversary_problem}
can be polynomially reduced to~\eqref{eq:follower_problem} in case of
discrete uncertainty by enumerating all scenarios and
solving~\eqref{eq:follower_problem} for each of them. For our setting,
in which~\eqref{eq:follower_problem} is solvable in polynomial time,
this immediately implies the following result:

\begin{theorem}\label{theo:discrete}
  Let $\U$ be a finite set, given explicitly as part of the
  input. Then the adversary's problem~\eqref{eq:adversary_problem} is
  solvable in polynomial time and the robust
  problem~\eqref{eq:robust_problem} is NP-easy.
\end{theorem}

When replacing the uncertainty set $\U$ in Theorem~\ref{theo:discrete}
by its convex hull, we obtain a significantly harder problem:

\begin{theorem}\label{theo:conv}
  Let~$\U$ be defined as the convex hull of a finite set of
  vectors, which are given explicitly as input. Then the robust
  counterpart~\eqref{eq:robust_problem} can be $\sigp 2$-hard and the
  adversary's problem~\eqref{eq:adversary_problem} NP-hard, even
  if~$X=\{0,1\}^p$.
\end{theorem}
\begin{proof}
  In the proof of Theorem~\ref{theorem:qsat}, we have seen that the
  linearized version of~\eqref{eq:qsat} is $\sigp 2$-hard for
  uncertainty sets of the form~$[-1,1]^n\times\{\bar c\}$, where~$\bar
  c$ is a vector collecting all certain entries of~$c$, corresponding
  to linearization variables and to artificial variables~$\bar y_i$ in
  the pessimistic case. It is easy to verify that the problem does not
  change when replacing~$[-1,1]^n$ by a superset, e.g., by the
  simplex~$\text{conv}\{-e,-e+2ne_1,\dots,-e+2ne_n\}$, where~$e$ is
  the all-ones vector and~$e_i$ is the $i$-th unit vector.  Since the
  latter is defined by polynomially many vertices, which can all be
  extended by~$\bar c$, we obtain the desired result
  for~\eqref{eq:robust_problem}. The NP-hardness
  of~\eqref{eq:adversary_problem} again follows directly from the
  $\sigp 2$-hardness of~\eqref{eq:robust_problem}.
\end{proof}

\section{Conclusion}

We have shown that bilevel optimization problems become significantly
har\-der when considering uncertain follower's objective functions in
a robust optimization approach. However, our results highlight a
fundamental difference between the interval uncertainty case and the
case of discrete uncertainty. Indeed, the construction in
Theorem~\ref{theorem:qsat} yields a problem
where~\eqref{eq:robust_problem} is $\sigp 2$-hard
and~\eqref{eq:adversary_problem} is NP-hard, while both problems are
at least one level easier in the setting of
Theorem~\ref{theo:discrete}. In this sense, interval uncertainty
renders bilevel optimization significantly harder than discrete
uncertainty, in contrast to classical single-level robust
optimization.  However, Theorem~\ref{theorem:qsat} and
Corollary~\ref{cor:discrete} show that~\eqref{eq:robust_problem} can
be one level harder than~\eqref{eq:basic_problem} for both interval
and discrete uncertainty.

Moreover, Theorem~\ref{theo:discrete} and Theorem~\ref{theo:conv}
together imply that the complexity of both the robust
counterpart~\eqref{eq:robust_problem} and the adversary's
problem~\eqref{eq:adversary_problem} may increase when replacing~$\U$
by its convex hull. This again is in contrast to the single-level
case, where it is well-known that replacing~$\U$ by its convex hull
essentially does not change the robust counterpart.

\bibliographystyle{plain}
\bibliography{literature_bilevel_robust}
  
\end{document}